\documentclass[12pt,reqno]{amsart}
\usepackage{amsmath, amsfonts, amssymb, amsthm, hyperref}
\usepackage{bm}
\allowdisplaybreaks[4]
\def\l{\left}
\def\r{\right}
\def\bg{\bigg}
\def\({\bg(}
\def\){\bg)}

\def\f{\frac}

\def\ls{\leqslant}
\def\gs{\geqslant}

\def\bi{\binom}

\def\eq{\equiv}

\def\da{\delta}

\def\Proof{\noindent{\it Proof}}

\def\Z{\mathbb Z}

\def\N{\mathbb N}

\def\1{{\bf 1}}

\def\pmod #1{\ ({\rm{mod}}\ #1)}

\def\<{\langle}
\def\>{\rangle}
\theoremstyle{plain}
\newtheorem{theorem}{Theorem}[section]
\newtheorem{lemma}{Lemma}
\newtheorem{corollary}{Corollary}

\theoremstyle{definition}
\newtheorem*{Ack}{Acknowledgment}

\theoremstyle{remark}
\newtheorem{remark}{Remark}

\makeatletter
\@namedef{subjclassname@2020}{%
  \textup{2020} Mathematics Subject Classification}
\makeatother
 \vspace{4mm}

\begin{document}
\hbox{Frontiers in Combinatorics and Number Theory 2 (2026), 55--62.}
\medskip

\title[Supercongruences for central trinomial coefficients]
{Supercongruences for central trinomial coefficients}
\author{Hao Pan}
\address {(Hao Pan) School of Applied Mathematics, Nanjing University of Finance and Economics,
 Nanjing 210046, People's Republic of China}
\email{haopan79@zoho.com}

\author{Zhi-Wei Sun}
\address {(Zhi-Wei Sun) School of Mathematics, Nanjing
University, Nanjing 210093, People's Republic of China}
\email{zwsun@nju.edu.cn}

\thanks{The paper was first posted to arXiv in 2020 as a preprint with the code arXiv:2012.05121. Both authors are supported by the National Natural Science Foundation of China (grants 12471312 and 12371004, respectively)}
\begin{abstract}
For each $n=0,1,2,\ldots$, the central trinomial coefficient $T_n$ is the coefficient of $x^n$
in the expansion of $(x^2+x+1)^n$. Let $p>3$ be a prime, and let $n$ be any positive integer.
In 2016, the second author conjectured that the quotient $(T_{pn}-T_n)/(pn)^2$ is always a $p$-adic integer. In this paper, we confirm this conjecture, and further prove that
$$\frac{T_{pn}-T_n}{(pn)^2}\equiv\frac{T_{n-1}}6\left(\frac p3\right)B_{p-2}\left(\frac13\right)\pmod p,$$
where $(\frac p3)$ is the Legendre symbol and $B_{p-2}(x)$ is the Bernoulli polynomial of degree $p-2$.
\end{abstract}
\maketitle

\section{Introduction}
\setcounter{lemma}{0}
\setcounter{theorem}{0}
\setcounter{equation}{0}
\setcounter{conjecture}{0}
\setcounter{remark}{0}
\setcounter{corollary}{0}

Let $b,c\in\Z$ and $n\in\N=\{0,1,2,\ldots\}$. As in \cite{S14a,S14b}, we denote by $T_n(b,c)$ the coefficient of $x^n$ in the expansion of $(x^2+bx+c)^n$. It is easy to see that
\begin{equation}\label{T_n(b,c)}T_n(b,c)=\sum_{k=0}^{\lfloor n/2\rfloor}\bi n{2k}\bi{2k}kb^{n-2k}c^k.
\end{equation}
When $b^2-4c\not=0$, it is known that
\begin{equation}\label{T-P}T_n(b,c)=(\sqrt{b^2-4c})^nP_n\l(\f b{\sqrt {b^2-4c}}\r),
\end{equation}
where
$$P_n(x):=\sum_{k=0}^n\bi nk\bi{n+k}k\l(\f{x-1}2\r)^k$$
is the Legendre polynomial of degree $n$.

 Those $T_n:=T_n(1,1)\ (n\in\N)$ are called {\it central trinomial coefficients}.
 Such numbers play important roles in enumerative combinatorics; for example,
 $T_n$ is the number of lattice paths from the point $(0, 0)$ to $(n, 0)$
with only allowed steps $(1, 1)$, $(1, -1)$ and $(1, 0)$ (cf. \cite[A002426]{Sl}).
For some known congruences involving central trinomial coefficients
or their generalizations, the reader may consult \cite{CS,CW,Liu,MS,S14a,S14b,S22}.

In 2016 the second author conjectured (cf. \cite[A277640]{S}) that for any prime $p>3$ and
$n\in\Z^+=\{1,2,3,\ldots\}$ we have
\begin{equation}\label{T_pn}\f{T_{pn}-T_n}{(pn)^2}\in\Z_p,\end{equation}
where $\Z_p$ denotes the ring of $p$-adic integers.
In this paper we confirm this conjecture and establish the following further result.

 \begin{theorem}\label{Th1.1}
Let $p>3$ be a prime and let $n$ be a positive integer. Then
\begin{equation}\label{T-cong}
\f{T_{pn}-T_n}{(pn)^2}\eq\f{T_{n-1}}6\l(\f p3\r)B_{p-2}\l(\f13\r)\pmod p,
\end{equation}
where $(\f p3)$ is the Legendre symbol and $B_{p-2}(x)$ is the Bernoulli polynomial of degree $p-2$.
\end{theorem}

This theorem has the following consequence.

\begin{corollary}\label{Cor1.1} For any prime $p>3$ and $a\in\Z^+$, we have
\begin{equation}\f{T_{p^a}-T_{p^{a-1}}}{p^{2a}}\eq\f16\l(\f p3\r)^aB_{p-2}\l(\f13\r)\pmod p.
\end{equation}
\end{corollary}

\begin{remark} By Corollary \ref{Cor1.1}, for any prime $p>3$ and $a\in\Z^+$ we have
$$T_{p^a}\eq T_{p^{a-1}}\eq\cdots\eq T_{p^0}=1\pmod{p^2}.$$
In 2016, the second author conjectured (cf. Comments in \cite[A002426]{Sl}) that for any integer $n>3$, the congruence
$T_n\eq 1\pmod{n^2}$ holds if and only if $n$ is prime.
This conjecture, if true,  provides an interesting characterization of primes via central trinomial coefficients.
\end{remark}

For $n=0,1,2,\ldots$, we define
\begin{equation}\label{a_n-def}a_n:=\sum_{k=0}^n\bi nk\bi{n-k}k\bi{n+k}k.
\end{equation}
Note that
\begin{equation}\label{a_n=}a_n=\sum_{k=0}^n\bi{n+k}{2k}\bi{n-k}k\bi{2k}k=\sum_{k=0}^n\bi {n+k}{3k}\bi{3k}k\bi{2k}k.
\end{equation}
For $|x|<1$, we have
\begin{align*}&\sum_{k=0}^\infty\bi{2k}k\bi{3k}k\f{x^{2k}}{(1-x)^{3k+1}}
\\=&\sum_{k=0}^\infty\bi{2k}k\bi{3k}kx^{2k}\sum_{l=0}^\infty\bi{3k+l}{3k}x^l
\\=&\sum_{n=0}^\infty x^n\sum_{k=0}^{\lfloor n/2\rfloor}\bi{2k}k\bi{3k}k\bi{n+k}{3k}=\sum_{n=0}^\infty a_nx^n.
\end{align*}
So the sequence $(a_n)_{n\gs0}$ coincides with \cite[A208425]{H}. Our second theorem
on this sequence confirms a conjecture of the second author (cf. Comments in \cite[A208425]{H})
formulated in 2016.

\begin{theorem} \label{Th1.2} Let $p>3$ be a prime and let $n$ be a positive integer. Then
\begin{equation}\label{a_n}\f{a_{pn}-a_n}{(pn)^3}\in\Z_p.\end{equation}
\end{theorem}

We are going to prove Theorem \ref{Th1.1} and Corollary \ref{Cor1.1}
in the next section. Theorem \ref{Th1.2} will be proved in Section 3.
Our proofs involve the cubic root
$$\omega:=\f{-1+\sqrt{-3}}2$$
of unity.

The reader may consult \cite{CD,OSS,S16} for some other similar known supercongruences,
and \cite[Conjecture 82]{S19} for some conjectural supercongruences.

\section{Proofs of Theorem \ref{Th1.1} and Corollary \ref{Cor1.1}}
\setcounter{lemma}{0}
\setcounter{theorem}{0}
\setcounter{equation}{0}
\setcounter{conjecture}{0}
\setcounter{remark}{0}
\setcounter{corollary}{0}

\begin{lemma}\label{Lem2.1} For any $n\in\N$ we have
\begin{equation}\label{T_n=} (-1)^nT_n=\sum_{k=0}^n\bi nk^2\omega^{2k-n}=\sum_{k=0}^n\bi nk^2\bar\omega^{2k-n}.
\end{equation}
\end{lemma}
\Proof. For $n\in\N$ it is known (cf. \cite[(3.134)]{G}) that
$$P_n(x)=\sum_{k=0}^n\bi nk^2\l(\f{x+1}2\r)^{n-k}\l(\f{x-1}2\r)^k.$$
Combining this with \eqref{T-P} we see that
\begin{align*}T_n=&\sqrt{-3}^nP_n\l(\f1{\sqrt{-3}}\r)=\sqrt{-3}^n\sum_{k=0}^n\bi nk^2\l(\f{1/\sqrt{-3}+1}2\r)^{n-k}\l(\f{1/\sqrt{-3}-1}2\r)^k
\\=&\sum_{k=0}^n\bi nk^2\l(\f{1+\sqrt{-3}}2\r)^{n-k}\l(\f{1-\sqrt{-3}}2\r)^k.
\end{align*}
Clearly,
$$\f{1+\sqrt{-3}}2=-\bar\omega=-\omega^2.$$
So, from the above we obtain
$$T_n=(-1)^n\sum_{k=0}^n\bi nk^2\omega^{2n-k}.$$
Taking conjugates of both sides, we get
$$T_n=(-1)^n\sum_{k=0}^n\bi nk^2\bar\omega^{2n-k}.$$
As $2n-k\eq 2k-n\pmod3$, we see that \eqref{T_n=} holds. \qed

\begin{lemma}\label{Lem2.2} Let $p>3$ be a prime. Then
\begin{equation}\label{3|k+p}\sum_{0<k<p\atop 3\mid k+p}\f1{k^2}
\eq-\f19\l(\f p3\r)B_{p-2}\l(\f13\r)\pmod p.
\end{equation}
\end{lemma}
\Proof. Note that
$$\sum_{0<k<p\atop 3\mid k+p}\f1{k^2}+\sum_{0<k<p\atop 3\mid k-p}\f1{k^2}+\sum_{0<k<p\atop 3\mid k}\f1{k^2}=\sum_{k=1}^{p-1}\f1{k^2}\eq0\pmod p$$
by Wolstenholme's congruence \cite{W}. As
$$\sum_{0<k<p\atop 3\mid k-p}\f1{k^2}\eq\sum_{0<k<p\atop 3\mid p-k}\f1{(p-k)^2}=\sum_{0<j<p\atop 3\mid j}\f1{j^2}\pmod{p},$$
by the above we obtain
\begin{equation}\label{2/9}\sum_{0<k<p\atop 3\mid k+p}\f1{k^2}\eq-2\sum_{0<k<p\atop 3\mid k}\f1{k^2}
=-\f 29\sum_{j=1}^{\lfloor (p-1)/3\rfloor}\f1{j^2}\pmod p.\end{equation}

By \cite[(9)]{Lehmer},
$$\sum_{j=1}^{\lfloor (p-1)/3\rfloor}\f1{j^2}\eq\sum_{j=1}^{\lfloor (p-1)/3\rfloor}j^{p-3}
\eq-\f1{p-2}B_{p-2}\l(\l\{\f p3\r\}\r)\pmod p,$$
where $\{p/3\}$ is the fractional part of $p/3$.
 Note that
 $$B_{p-2}\l(\l\{\f p3\r\}\r)=\l(\f p3\r)B_{p-2}\l(\f 13\r)$$
 since
$$B_{p-2}\l(\f 23\r)=(-1)^{p-2}B_{p-2}\l(\f13\r)$$
(cf. \cite[p.\,248]{IR}). Therefore
$$\sum_{j=1}^{\lfloor(p-1)/3\rfloor}\f1{j^2}\eq\f12\l(\f p3\r)B_{p-2}\l(\f 13\r)\pmod p$$
and hence by \eqref{2/9} we have
$$\sum_{0<k<p\atop 3\mid k+p}\f1{k^2}\eq-\f29\sum_{j=1}^{\lfloor(p-1)/3\rfloor}\f1{j^2}
\eq-\f19\l(\f p3\r)B_{p-2}\l(\f 13\r)\pmod p.$$
This concludes the proof. \qed

\medskip
\noindent
{\it Proof of Theorem \ref{Th1.1}}. Let $k\in\{0,\ldots,pn\}$ with $p\nmid k$. Then
\begin{align*}\bi{pn}k=&\f{pn}k\prod_{0<j<k}\f{pn-j}j
\\=&\f{pn}k\prod_{0<i\ls\lfloor(k-1)/p\rfloor}\f{pn-pi}{pi}\times\prod_{0<j<k\atop p\nmid j}\l(\f{pn}j-1\r)
\\\eq&\f{pn}k\bi{n-1}{\lfloor(k-1)/p\rfloor}(-1)^{\da_k}\pmod {p^{2+2\nu_p(n)}}
\end{align*}
where $\da_k:=|\{0<j<k:\ p\nmid j\}$ and $\nu_p(n):=\max\{a\in\N:\ p^a\mid n\}$.
It follows that
\begin{equation}\label{pnk}\sum_{0\ls k\ls pn\atop p\nmid k}\bi{pn}k^2\omega^{2k-pn}
\eq\sum_{0\ls k\ls pn\atop p\nmid k}\f{p^2n^2}{k^2}\bi{n-1}{\lfloor(k-1)/p\rfloor}^2\omega^{2k-pn}
\pmod{p^{3+2\nu_p(n)}}.
\end{equation}
By Lemma \ref{Lem2.1},
$$(-1)^{n-1}T_{n-1}=\sum_{k=0}^{n-1}\omega^{2k-n+1}\bi{n-1}k^2
=\sum_{k=0}^{n-1}\bar\omega^{2k-n+1}\bi{n-1}k^2.$$
As $p\eq\pm1\pmod 3$, we get
$$\sum_{k=0}^{n-1}\omega^{p(2k-n+1)}\bi{n-1}k^2=(-1)^{n-1}T_{n-1}.$$
Therefore
\begin{align*}&\sum_{0\ls k\ls pn\atop p\nmid k}\f{\omega^{2k-pn}}{k^2}\bi{n-1}{\lfloor(k-1)/p\rfloor}^2
\\=&\omega^{-pn}\sum_{k=0}^{n-1}\sum_{j=1}^{p-1}\f{\omega^{2(pk+j)}}{(pk+j)^2}\bi{n-1}k^2
\eq\omega^{-pn}\sum_{k=0}^{n-1}\omega^{2pk}\bi{n-1}k^2\sum_{j=1}^{p-1}\f{\omega^{2j}}{j^2}
\\\eq&\omega^{-p}\sum_{k=0}^{n-1}\omega^{p(2k-n+1)}\bi{n-1}k^2\sum_{j=1}^{p-1}\f{\omega^{2j}}{j^2}
=\omega^{-p}(-1)^{n-1}T_{n-1}\sum_{j=1}^{p-1}\f{\omega^{2j}}{j^2}\pmod{p}.
\end{align*}
In light of this and the equality $\omega+\bar\omega=-1$, we get
\begin{align*}&\sum_{0\ls k\ls pn\atop p\nmid k}\f{\omega^{2k-pn}+\bar \omega^{2k-pn}}{k^2}\bi{n-1}{\lfloor(k-1)/p\rfloor}^2
\\\eq&(-1)^{n-1}T_{n-1}\sum_{j=1}^{p-1}\f{\omega^{2j-p}+\bar\omega^{2j-p}}{j^2}
\eq (-1)^{n-1}T_{n-1}\(\sum_{0<j<p\atop 3\mid 2j-p}\f2{j^2}-\sum_{0<j<p\atop 3\nmid 2j-p}\f1{j^2}\)
\\\eq&(-1)^{n-1}T_{n-1}\(\sum_{0<j<p\atop 3\mid 2j-p}\f3{j^2}-\sum_{j=1}^{p-1}\f1{j^2}\)
\eq3(-1)^{n-1}T_{n-1}\sum_{0<j<p\atop 3\mid j+p}\f1{j^2}\pmod p.
\end{align*}
Combining this with \eqref{pnk} and Lemma \ref{Lem2.2}, we obtain
\begin{equation}\label{nmid}\sum_{0\ls k\ls pn\atop p\nmid k}\bi{pn}k^2\f{\omega^{2k-pn}+\bar \omega^{2k-pn}}{2}
\eq(-1)^n\f{T_{n-1}}6\l(\f p3\r)(pn)^2B_{p-2}\l(\f13\r)\pmod {p^{3+2\nu_p(n)}}.
\end{equation}

In view of Lemma \ref{Lem2.1} and \eqref{nmid}, we have
\begin{align*}
(-1)^{pn}T_{pn}-(-1)^nT_n=&\sum_{k=0}^{pn}\bi{pn}k^2\f{\omega^{2k-pn}+\bar\omega^{2k-pn}}2-\sum_{k=0}^n\bi nk^2
\f{\omega^{2k-n}+\bar\omega^{2k-n}}2
\\=&\sum_{k=0}^n\bi{pn}{pk}^2\f{\omega^{p(2k-n)}+\bar\omega^{p(2k-n)}}2-\sum_{k=0}^n\bi nk^2\f{\omega^{2k-n}+\bar\omega^{2k-n}}2
\\&+\sum_{0\ls k\ls pn\atop p\nmid k}\bi{pn}k^2\f{\omega^{2k-pn}+\bar \omega^{2k-pn}}{2}
\\\eq&\sum_{k=0}^n\(\bi{pn}{pk}^2-\bi nk^2\)\f{\omega^{2k-n}+\bar\omega^{2k-n}}2
\\&+(-1)^n\f{T_{n-1}}6\l(\f p3\r)(pn)^2B_{p-2}\l(\f13\r)
\pmod {p^{3+2\nu_p(n)}}
\end{align*}
since $p\eq\pm1\pmod3$.
For each $k\in\{0,\ldots,n\}$, it is known (see, e.g., \cite{RZ}) that
\begin{equation}\label{pn-pk}\bi{pn}{pk}=\bi nk(1+p^3nk(n-k)q_k)
\end{equation} for some $q_k\in\Z_p$, hence
$$\bi{pn}{pk}^2-\bi nk^2=2p^3\bi nk^2nk(n-k)q_k+p^6\l(\bi nknk(n-k)q_k\r)^2
\eq0\pmod{p^{3+3\nu_p(n)}}$$
since $\bi nk k=n\bi{n-1}{n-k}$ and $\bi nk(n-k)=n\bi{n-1}k$.
So we finally get
$$T_{pn}-T_n\eq\f{T_{n-1}}6\l(\f p3\r)(pn)^2B_{p-2}\l(\f13\r)
\pmod {p^{3+2\nu_p(n)}}$$
which is equivalent to \eqref{T-cong}. \qed

\medskip
\noindent{\it Proof of Corollary 1.1}. Applying Theorem \ref{Th1.1} with $n=p^{a-1}$, we get
$$\f{T_{p^a}-T_{p^{a-1}}}{p^{2a}}\eq\f{T_{p^{a-1}-1}}6\l(\f p3\r)B_{p-2}\l(\f13\r)\pmod p.$$
So it suffices to prove
\begin{equation}\label{a-1}T_{p^{a-1}-1}\eq\l(\f p3\r)^{a-1}\pmod p.
\end{equation}

Clearly \eqref{a-1} holds for $a=1$. Below we assume $a>1$. For each $k=1,\ldots,(p^{a-1}-1)/2$, we have
$$\bi{p^{a-1}-1}{2k}=\prod_{j=1}^{2k}\l(\f{p^{a-1}}j-1\r)\eq(-1)^{2k}=1\pmod p$$
and
$$\f{\bi{(p^{a-1}-1)/2}k}{\bi{-1/2}k}=\prod_{0\ls j<k}\l(1-\f{p^{a-1}}{2j+1}\r)\eq1\pmod p.$$
Thus
\begin{align*}T_{p^{a-1}-1}=&\sum_{k=0}^{(p^{a-1}-1)/2}\bi{p^{a-1}-1}{2k}\bi{2k}k
=\sum_{k=0}^{(p^{a-1}-1)/2}\bi{p^{a-1}-1}{2k}\bi{-1/2}k(-4)^k
\\\eq&\sum_{k=0}^{(p^{a-1}-1)/2}\bi{(p^{a-1}-1)/2}{k}(-4)^k
=(1-4)^{(p^{a-1}-1)/2}=(-3)^{\f{p-1}2\sum_{0\ls r<a-1}p^r}
\\\eq&\l(\f{-3}p\r)^{a-1}=\l(\f p3\r)^{a-1}\pmod p
\end{align*}
with the aid of the theory of quadratic residues.
This proves the desired \eqref{a-1}.

In view of the above, we have completed the proof of Corollary \ref{Cor1.1}. \qed

\section{Proof of Theorem \ref{Th1.2}}
\setcounter{lemma}{0}
\setcounter{theorem}{0}
\setcounter{equation}{0}
\setcounter{conjecture}{0}
\setcounter{remark}{0}
\setcounter{corollary}{0}

\begin{lemma} For any $n\in\N$ we have
\begin{equation}\label{3.1}a_n=\sum_{k=0}^n\bi nk^3(-\omega)^{2k-n}
=\sum_{k=0}^n\bi nk^3(-\bar\omega)^{2k-n}.
\end{equation}
\end{lemma}
\Proof. Recall MacMahon's identity \cite[(6.7)]{G}
$$\sum_{k=0}^n\bi nk^3x^k=(1+x)^n\sum_{k=0}^{\lfloor n/2\rfloor}\bi{n+k}{3k}\bi{2k}k\bi{3k}k\l(\f x{(1+x)^2}\r)^k.$$
Putting $x=\omega^2$ and noting $(x+1)^2=x$, we get
$$\sum_{k=0}^n\bi nk^3\omega^{2k}=(-\omega)^n\sum_{k=0}^n\bi {n+k}{3k}\bi{2k}k\bi{3k}k$$
and hence
$$\sum_{k=0}^n\bi nk^3(-\omega)^{2k-n}=a_n$$
by \eqref{a_n=}.
Taking conjugates of both sides of the last equality, we obtain
$$\sum_{k=0}^n\bi nk^3(-\bar\omega)^{2k-n}=a_n.$$
So \eqref{3.1} is valid. \qed
\medskip

\noindent{\it Proof of Theorem 1.2}. For $k\in\{0,\ldots,n\}$, we have \eqref{pn-pk} for some $q_k\in\Z_p$, hence
$$\bi{pn}{pk}^3=\bi nk^3(1+p^3nk(n-k)q_k)^3\eq\bi nk^3\pmod{p^{3+3\nu_p(n)}}$$
since $k\bi nk=n\bi{n-1}{n-k}$ and $(n-k)\bi nk=n\bi{n-1}k$.
Thus
\begin{align*}a_{pn}=&\sum_{k=0}^n\bi {pn}{pk}^3(-\omega)^{2pk-pn}+\sum_{0<k<pn\atop p\nmid k}\f{(pn)^3}{k^3}\bi{pn-1}{k-1}^3(-\omega)^{2k-pn}
\\\eq&\sum_{k=0}^n\bi nk^3(-\omega)^{(\f p3)(2k-n)}(-\omega)^{(p-(\f p3))(2k-n)}
\\=&\sum_{k=0}^n\bi nk^3(-\omega^{(\f p3)})^{2k-n}=a_n\pmod{p^{3+3\nu_p(n)}}
\end{align*}
and hence \eqref{a_n} follows. \qed

\Ack The authors are indebted to the referee for helpful comments.

\medskip


\end{document}